\documentclass[a4]{article}
\usepackage{fullpage}
\usepackage{amsfonts}
\usepackage{amsmath}
\usepackage{amssymb}
\usepackage{graphicx}
\usepackage{xcolor}
\usepackage[mathscr]{eucal}

\title{On $SL(2,\mathbb{R})$-cocycles over irrational rotations with secondary collisions}
\author{Alexey V. Ivanov}
\date{}
\begin{document}
\renewcommand{\theequation}{\arabic{section}.\arabic{equation}}
\maketitle

\begin{abstract}
We consider a skew product $F_{A} = (\sigma_{\omega}, A)$ over irrational rotation $\sigma_{\omega}(x) = x + \omega$  of a circle $\mathbb{T}^{1}$. It is supposed that the transformation $A: \mathbb{T}^{1} \to SL(2, \mathbb{R})$ being a $C^{1}$-map has the form $A(x) = R(\varphi(x)) Z(\lambda(x))$, where $R(\varphi)$ is a rotation in $\mathbb{R}^{2}$ over the angle $\varphi$ and $Z(\lambda)= diag\{\lambda, \lambda^{-1}\}$ is a diagonal matrix. Assuming that $\lambda(x) \ge \lambda_{0} > 1$ with a sufficiently large constant $\lambda_{0}$ and the function $\varphi$ be such that $\cos \varphi(x)$ possesses only simple zeroes, we study hyperbolic properties of the cocycle generated by $F_{A}$. We apply the critical set method to show that, under some additional requirements on the derivative of the function $\varphi$, the secondary collisions compensate weakening of the hyperbolicity due to primary collisions and the cocycle generated by $F_{A}$ becomes hyperbolic in contrary to the case when secondary collisions can be partially eliminated.
\end{abstract}

Keywords: linear cocycle, hyperbolicity, Lyapunov exponent, critical set

MSC 2010: 37C55, 37D25, 37B55, 37C60

\section{Introduction}

One of the fundamental problems in the theory of smooth dynamics is to establish for a given dynamical system which hyperbolic properties it possesses. In particular, it is  important to determine whether the system is uniformly hyperbolic or not, and in the latter case, to prove (or disprove) the positiveness of its Lyapunov exponents. Often, instead of one particular system, a family of dynamical systems is considered. Thus, it naturally arises a necessity to answer the questions mentioned above with respect to parameter values of the family. There is wide literature on this subject (see e.g. \cite{BarPes}, \cite{BDV} and references therein). One may note that the difficulty of the problem increases together with growth of dimension of considered dynamical systems. Due to this fact, the case of one-dimensional discrete systems is much more explored in comparison even with two-dimensional case. An intermediate position is occupied by skew-products over one-dimensional cascades: they keep many features of multidimensional systems, but the governing dynamics is much easier. They have been studied for several decades and from different points of view. We refer the reader to the following (of course, far from to be complete) list of papers: \cite{Avi}, \cite{AvKri}, \cite{BourJit}, \cite{Eli}, \cite{Her}, \cite{JohnMos}, \cite{SorSpe}. However, it has to be noted that there is a lack of direct constructive methods, which allow to solve the mentioned problem for a particular dynamical system or a family of systems (see e.g. \cite{AvBo_01}, \cite{BuFe}, \cite{Her}, \cite{LSY}). 

In this paper we continue a study of skew-products over irrational rotation started in \cite{Iva21}. Let
\begin{equation}\label{eq1}
F_{A}: \mathbb{T}^{1}\times \mathbb{R}^{2} \to \mathbb{T}^{1}\times \mathbb{R}^{2}
\end{equation}
be a skew-product map defined by 
$$
(x, v)\mapsto (\sigma_{\omega}(x), A(x)v),\;
(x, v) \in \mathbb{T}^{1}\times \mathbb{R}^{2}
$$
where $\sigma_{\omega}(x) = x + \omega$ is a rotation of a circle $\mathbb{T}^{1}=\mathbb{R}/\mathbb{Z}$ with an irrational rotation number $\omega$ and 
$$
A: \mathbb{T}^{1} \to SL(2, \mathbb{R})
$$ 
is a $C^{1}$-function.  

Interest to such skew-products is determined not only by the fact that they can be considered as a bridge between one- and two-dimensional cascades, but also their direct applications in physics. One may associate to (\ref{eq1}) the following difference equation
\begin{equation}\label{eqSys}
\psi(y+\omega) = \mathcal{A}(y) \psi(y),\; y\in \mathbb{R}.
\end{equation}
Here $\mathcal{A} = A\circ \pi_{st}$ is a 1-periodic matrix-valued function, $\pi_{st}:\mathbb{R}\to \mathbb{T}^{1}=\mathbb{R}/\mathbb{Z}$ is the quotient map and $\psi = (\psi_{1}, \psi_{2})^{tr}$ is an unknown vector-function.
Eliminating the second component $\psi_{2}$ leads to a second-order difference equation for $\psi_{1}$:
\begin{equation}\label{eqDif}
m(y)\psi_{1}(y+\omega) + n(y)\psi_{1}(y) + p(y)\psi_{1}(y-\omega) = 0.
\end{equation}
Here $m, n, p$ are known $1$-periodic real-valued functions expressed in terms of entries of $\mathcal{A}$.

Such difference equations have many applications. In particular, they appear in spectral theory of the Schr\"odinger operator on $l_{2}(\mathbb{Z})$ and in the stability problem for the Hill's equation with quasiperiodic potentials \cite{Avi}, \cite{ABD}, \cite{Iva17}. Another application comes from the electromagnetic-wave diffraction in a wedge-shaped region. Indeed, the Sommerfeld-Malyuzhinets representation for the electic field leads to a system of linear difference equations for two coupled spectral functions (see e.g. \cite{Lya}). Eliminating one of them gives a second-order difference equation of type (\ref{eqDif}) for the remaining spectral function.

It is remarkable fact that the property of equation (\ref{eqDif}) to possess a solution in one or another functional space correlates with dynamical properties of the corresponding skew-product 
\cite{Puig}, \cite{Sack}.
For example, the property of exponential dichotomy for (\ref{eqSys}) is equivalent to the uniform hyperbolicity for the skew-product \cite{AvBo}. 
 
In the present paper we consider a skew-product (\ref{eq1}) satisfying the following assumptions. Namely, we suppose that the transformation $A$ can be represented as
\begin{equation}\label{eqA}
A(x) = R(\varphi(x))\cdot Z(\lambda(x)), 
\end{equation}
where
$$
R(\varphi) = \left(
\begin{array}{cc}
\cos \varphi & \sin \varphi\\
-\sin \varphi & \cos \varphi
\end{array}
\right), 
\,
Z(\lambda) = \left(
\begin{array}{cc}
\lambda & 0\\
0 & \lambda^{-1}
\end{array}
\right)
$$ 
with some $C^{1}$-functions $\varphi: \mathbb{T}^{1}\to 2\pi\mathbb{T}^{1}$, $\lambda: \mathbb{T}^{1}\to \mathbb{R}$ such that
\begin{align*}
(H_{1})\quad & \left\{x\in \mathbb{T}^{1}:\; \cos(\varphi(x)) = 0\right\} = \bigcup\limits_{j=0}^{N} \{c_{j}\},\\
(H_{2})\quad & \forall\,\,  j=0,\ldots, N,\,\, C_{1}^{(j)}\varepsilon^{-1}\le \vert\varphi'(x)\vert \le C_{2}^{(j)}\varepsilon^{-1},\,\, \forall x\in U_{\varepsilon}(c_{j}),\,\, \varepsilon\ll 1;\,\,
Var(\phi,  U_{\varepsilon}(c_{j})) \sim O(1),\\
(H_{3})\quad & \bigl\vert \cos(\varphi(x))\bigr\vert \ge C_{3}\quad
\forall x\in \mathbb{T}^{1}\setminus \bigcup\limits_{j=0}^{N} U_{\varepsilon}(c_{j}),\\
(H_{4})\quad & ind(\varphi) = 0,\\
(H_{5})\quad & \lambda(x)\ge \lambda_{0} > 1 \quad \forall x\in \mathbb{T}^{1}.
\end{align*}
Here and after $C_{k}$ denotes a positive constant, $U_{\varepsilon}(x)$ is the $\varepsilon$-neighbourhood of a point $x\in \mathbb{T}^{1}$, $Var(\phi,  U_{\varepsilon}(c_{j}))$ is the variation of the function $\varphi$ in the neighborhood $U_{\varepsilon}(c_{j})$ and $ind(\varphi)$ stands for the index of a closed curve $\varphi(\mathbb{T}^{1})$.
Additionally, we assume that functions $\varphi$ and $\lambda$ depend smoothly on a real parameter $t\in [a, b]\subset \mathbb{R}$ such that
\begin{align*}
&(H_{6})\quad \left\vert\frac{\textrm{d}\rho(c_{j}(t), c_{k}(t))}{\textrm{d} t}\right\vert > C_{4} > 0,\,\, \forall t\in [a, b],\,\, j\neq k,
\end{align*}
where $\rho$ denotes the standard distance in $\mathbb{T}^{1}$.

To the skew product (\ref{eq1}) one may assign a cocycle $M(x,n)$ defined as
\begin{align*}
&M(x, n)=A(\sigma_{\omega}^{n-1}(x))\ldots A(x), \,\, n>0;\\ 
&M(x, n)=\left[A(\sigma_{\omega}^{-n}(x))\ldots A(\sigma_{\omega}^{-1}(x))\right]^{-1}, \,\, n<0;\\
&M(x, 0)=I.
\end{align*}

In this paper we study how the property of hyperbolicity is related to the cocycle parameters.
\newtheorem{defs}{Definition}
\begin{defs}
We say that a cocycle $M$ is uniformly hyperbolic (UH) if there exist continuous maps $E^{u,s}: \mathbb{T}^{1} \to Gr(2,1)$ and positive constants $C, \Lambda$  such that the subspaces 
$E^{u,s}(x)$ are invariant with respect to the map (\ref{eq1})  
(i.e. $E^{u,s}(\sigma_{\omega}(x))=A(x)E^{u,s}(x)$) and $\forall x\in \mathbb{T}^{1}$, $n\ge 0$
\begin{align}
\left\Vert M(x, -n)\vert_{E^{u}(x)}\right\Vert \le C{\rm e}^{-\Lambda n}, \notag\\
\left\Vert M(x, n)\vert_{E^{s}(x)}\right\Vert \le C{\rm e}^{-\Lambda n}.\notag
\end{align}
\end{defs}
Here $Gr(2,1)$ stands for the set of $1$-dimensional subspaces of $\mathbb{R}^{2}$. We note that the Oseledets theorem guarantees the existence of such invariant subspaces for a.e. $x\in \mathbb{T}^{1}$, but, in general, the maps $E^{u,s}$ are only measurable.

For our purpose it is more convenient to use an alternative version of this definition (see e.g. \cite{AvBo}). Namely, the cocycle $M$ is said to be UH if there exist positive constants $C, \Lambda_{0}$ such that $\forall x\in \mathbb{T}^{1}$ and $n\ge 0$
\begin{equation}\label{eqUH}
\Vert M(x,n)\Vert \ge C {\rm e}^{\Lambda_{0} n}.
\end{equation}

It has to be noted that due to Kingman's subadditive ergodic theorem for a.e. $x\in \mathbb{T}^{1}$ there exists the Lyapunov exponent 
$$\Lambda(x) = \lim\limits_{n\to \infty}\frac{1}{n}\log\Vert  M(x,n)\Vert.$$
Moreover, since $\omega$ is assumed to be irrational, the rotation $\sigma_{\omega}$ is ergodic and
$$\Lambda(x) = \tilde\Lambda_{0}\quad a.e.,$$
where $\tilde\Lambda_{0}$ is the integrated Lyapunov exponent 
$$\tilde\Lambda_{0} = \int \Lambda(x){\rm d}x = \lim\limits_{n\to \infty}\int \frac{1}{n}\log\Vert  M(x,n)\Vert{\rm d}x.$$
Hence, the UH implies positiveness of $\tilde \Lambda_{0}$, but the opposite is not true (see e.g. \cite{Her}). One calls a cocycle $M$ non-uniformly hyperbolic (NUH) if $\tilde \Lambda_{0} > 0$, but $M$ is not UH. Uniformly hyperbolic cocycles constitute an open subset in the set of all cocycles. On the other hand, if a cocycle is NUH, its Lyapunov exponent can be made equal to zero by arbitrary small $C^{0}$-perturbation. Thus, even the problem of finding an example of NUH with positive Lyapunov exponent is not trivial. In \cite{Her} M. Herman provided the first example of such cocycle. The constructed cocycle corresponds to a skew-product of type (\ref{eqA}) such that function $\lambda$ is constant $\lambda(x) = \lambda_{0}>1$ and $ind(\varphi)=1$. The latter condition guarantees a topologic obstacle for uniform hyperbolicity. 

{\bf Remark 1} One may remark (see e.g. \cite{BenCar}, \cite{LSY}) that the problem admits the projectivization in the following sense. Consider the standard covering $(p, \mathbb{T}^{1}, \mathbb{R}P^{1})$ of the real projective line, where the projection $p$ identifies the diametrically opposite points of the circle. For any continuous function $\hat\varphi: \mathbb{T}^{1}\to \mathbb{T}^{1}$ it generates a continuous function $\hat\varphi_{p}=p\circ \varphi:\mathbb{T}^{1}\to \mathbb{R}P^{1}$. On the other hand, for a given continuous function 
$\hat\varphi_{p}: \mathbb{T}^{1}\to \mathbb{R}P^{1}$ its lift $\tilde\varphi$ is not necessary continuous. Moreover, there exists exactly two continuous lifts, which will be denoted by 
$\hat\varphi_{k}, k=1,2$. Let $\tilde\varphi_{k}, k=1,2$ be two arbitrary (not necessary continuous) lifts of $\hat\varphi_{p}$. It follows from Definition 1 that the cocycles, defined by (\ref{eq1}) and corresponding to these lifts, are either both UH or both not UH. Particularly, if they are both UH, the stable (unstable) subspaces $E^{s,(u)}_{k}(x)$ coincide $E^{s,(u)}_{1}(x)=E^{s,(u)}_{2}(x)$. This remark enables us to consider the function $\varphi$ from (\ref{eqA}) as a lift of some continuously differentiable function $\varphi_{p}: \mathbb{T}^{1}\to \mathbb{R}P^{1}$. 

In the present work we use an approach developed in \cite{Jak}, \cite{BenCar}, \cite{LSY}, \cite{Laz}. Devoted to different objects these papers have a common framework. The idea suggested initially in \cite{Jak} can be described as follows. Consider a family of skew-product maps $F_{A,t} = (f, A_{t})$ dependent on a parameter $t$ and defined (similarly to (\ref{eq1})) on a vector bundle $V$ over a base $B$. Properties  of the fiber transformation $A_{t}$ may vary with respect to a point of the base. Selecting those points of the base which correspond to violation of some specific property (e.g. hyperbolicity) of $A_{t}$, we construct the so-called critical set $\mathcal{C}_{0}$. Taking a small neighbourhood of $\mathcal{C}_{0}$, one needs to study the dynamics of this set under the map $f$.  However, due to dependence on the parameter $t$, interactions between different parts of $\mathcal{C}_{0}$ may have degeneracies. Detuning the parameter $t$, we exclude such degeneracies and put the interactions  in a general position. Finally, using properties of the base map $f$ (e.g. ergodicity) and non-degeneracy of interactions of the critical set, one may extract an additional information on the whole system.

The paper is organized as follows. In Section 2 we define the critical set and  introduce a notion of dominance for primary and secondary collisions. We also formulate necessary conditions which guarantee the dominance of primary collisions and, as a consequence, non-uniform hyperbolicity of the cocycle. In section 3 the case, which corresponds to dominance of secondary collisions for the simplest critical set consisting of two points, is considered. We perform asymptotic analysis as $\lambda_{0}\to +\infty$, $\varepsilon\to 0$ to describe effect of interaction between small neighborhoods of two critical points on the cocycle. Resonance conditions on the rotation number and parameter $t$ are presented. Under these conditions we prove the uniform hyperbolicity for the cocycle (\ref{eqA}).

\section{Dynamics of the critical set}

Note that, by definition, the cocycle $M$ corresponding to (\ref{eq1}) is a product of matrices $A_{k}$ such that $\Vert A_{k}\Vert\ge \lambda_{0}$ are sufficiently large.
However, the product may not admit estimate (\ref{eqUH}). The obstacle to this fact is the presence of the critical set, which can be defined for the skew product (\ref{eqA}) as
\begin{equation}\label{eqCoc}
\mathcal{C}_{0} = \bigcup\limits_{j=0}^{N}\{c_{j}\}.
\end{equation}
To understand a role of the critical set, we formulate the following technical lemma (see e.g. \cite{Iva21})
\newtheorem{lems}{Lemma}
\begin{lems}
For any fixed $\phi\in [0, \pi)$ and $\lambda_{1}, \lambda_{2}$ such that $\lambda_{k} > 1, k=1,2$ the following representation holds true
\begin{equation*}
Z(\lambda_{2}) R(\phi) Z(\lambda_{1}) = R(\psi - \chi) Z(\mu) R(\chi),
\end{equation*}
where $\psi\in [0, \pi)$, $\chi\in [-\pi/4, \pi/4]$, $\mu\ge 1$
\begin{align}
\nonumber
&\mu = \frac{a}{2}\left(1+ c + \left(\left(1+ c\right)^{2}-4 a^{-2}\right)^{1/2}\right),\\
\nonumber
&a=\frac{\lambda_{1}}{\lambda_{2}}\cdot \frac{\lambda_{2}^{2}\cos^{2}\phi + \beta\sin^{2}\phi}{\left(\cos^{2}\phi + \beta^{2}\sin^{2}{\phi}\right)^{1/2}},\\
\label{eq_Lem1}
&b=\left(\frac{\lambda_{2}}{\lambda_{1}}\right)^{2}\frac{1-\lambda_{2}^{-4}}{1+\lambda_{1}^{-2}\lambda_{2}^{-2}}\cdot
     \frac{\sin\phi \cos\phi}{\lambda_{2}^{2}\cos^{2}\phi + \beta\sin^{2}\phi},\\
\nonumber
&c=b^{2}+a^{-2},\quad \beta =  \frac{\lambda_{1}^{-2} + \lambda_{2}^{-2}}{1 + \lambda_{1}^{-2}\lambda_{2}^{-2}},\quad \tan\psi = \beta\tan\phi,\\ 
\nonumber
&\tan \chi = -\frac{\sqrt{2}b(1-c)^{-1}}
{\sqrt{1 + 2b^{2}(1-c)^{-2} + \sqrt{1+4b^{2}(1-c)^{-2}}}}.
\end{align}
\end{lems}

One may note that the norm of a product $P_{1} = Z(\lambda_{2}) R(\phi) Z(\lambda_{1})$ equals to $\mu$ and, by Lemma 1, 
\begin{equation*}
\mu = 1 \Leftrightarrow a = 1,\;\; b = 0.
\end{equation*}
Besides, if $\lambda_{1}, \lambda_{2}$ are sufficiently large, then we have the following implications:
\begin{align*}
&1.\quad \cos\phi \nsim 0\quad\Rightarrow\quad \mu\sim a\sim \lambda_{1}\lambda_{2}\bigl\vert\cos\phi\bigr\vert,\\
&2.\quad \cos\phi \sim 0\quad\Rightarrow\quad \mu\sim a\sim \frac{\lambda_{1}}{\lambda_{2}},\\
&3.\quad \phi = \frac{\pi}{2}\, \textrm{mod}\, \pi\quad \Rightarrow\quad \psi= \frac{\pi}{2}\, \textrm{mod}\, \pi,\; \chi = 0.
\end{align*}
It means that, if $\vert\cos \phi\vert$ is bounded away from zero, the norm of $P_{1}$ increases with respect to $\lambda_{1}$, since it is multiplied by a large quantity proportional to $\lambda_{2}\bigl\vert\cos\phi\bigr\vert$. On the other hand, if 
$\vert\cos \phi\vert$ is sufficiently small, the norm of $P_{1}$ decreases with respect to $\lambda_{1}$, as it is divided by a quantity proportional to $\lambda_{2}$. Moreover, in the latter case the angle $\psi$ satisfies $\cos\psi \sim 0$.  Applying this lemma to the cocycle $M$, we conclude that, each time a trajectory of a point $x\in \mathbb{T}^{1}$ under the rotation $\sigma_{\omega}$ falls into a small neighbourhood of the set $\mathcal{C}_{0}$, the growth of the cocycle norm changes its behaviour from increasing to decreasing or vice versa. This behaviour continues until the next visit of the trajectory to a small neighborhood of $\mathcal{C}_{0}$. Since $\sigma_{\omega}$ is ergodic, every point $x\in\mathbb{T}^{1}$  reaches a small neighbourhood of $\mathcal{C}_{0}$ after some iterations. Thus, the hyperbolic properties of $M(x,n)$ strongly depend on the dynamics of the critical set itself. 

To describe the dynamics of the critical set, we introduce for a fixed sufficiently small $\delta > 0$ a $\delta$-neighbourhood of $\mathcal{C}_{0}$ denoted by
\begin{equation*}
U_{\delta} = \bigcup\limits_{j=0}^{N}I_{j}(\delta),
\end{equation*}
where $I_{j}(\delta) = \{x\in \mathbb{T}^{1}: \rho(x, c_{j})<\delta\}$.

Then, due to ergodicity of $\sigma_{\omega}$, a trajectory of any point $c_{j}$ enters each interval $I_{k}(\delta)$ infinitely many times. It enables us to introduce notions of collision and time of collision.
\begin{defs}
Let $\tau_{j,j'}$ be the minimum of integer $k > 0$ such that
$$
\sigma_{\omega}^{k}(c_{j})\cap I_{j'} \neq \emptyset.
$$
We say that the points $c_{j}$ and $c_{j'}$ collide with accuracy $\delta$ at the time $\tau_{j,j'}$ and call such event a collision and $\tau_{j,j'}$ the time of collision.  
A collision is called primary if $j=j'$ and secondary if $j\neq j'$. 
\end{defs}
There is essential difference in behaviour of the primary and secondary collisions with respect to the parameter $t$. First, we note that the times of primary collisions $\tau_{j,j}$ do not depend on $j$  and we may denote them by $\tau_{0}$. It is a characteristic of the rotation number $\omega$ and the parameter $\delta$ only. On the other hand, the assumption ($H_{6}$) implies that relative positions of the points $c_{j}$ vary with respect to $t$ and, hence, the times of secondary collisions depend on the parameter $t$ in a non-trivial way.  

Following \cite{Iva21} we give a definition.
\begin{defs}
For a fixed $\delta>0$ we say that primary collisions dominate if for all pairs $(j, j')$ such that $j\neq j'$
$$
\tau_{j,j'}(\delta) > \tau_{0}(\delta).
$$
Otherwise, we say that secondary collisions dominate.
\end{defs}

In \cite{Iva21} we investigated the problem of elimination of the secondary collisions by detuning the parameters of the cocycle. Due to ergodicity of $\sigma_{\omega}$, one cannot eliminate 
these collisions completely. However, if the rotation number $\omega$ satisfies two number-theoretical conditions, it is possible to achieve domination of primary collision for some decreasing sequence $\{\delta_{k}\}$. To formulate these conditions, denote by $p_{n}/q_{n}$ the best rational approximation of order $n$ to $\omega$. We assume that $\omega$ satisfies the Brjuno's condition with a constant $C_{B}$ (see \cite{Brj}, \cite{Russ}), i.e.
\begin{equation}\label{eqCondB}
\sum\limits_{n=1}^{\infty}\frac{\log (2 q_{n+1})}{q_{n}} = C_{B} < \infty.
\end{equation}
The second condition can be formulated as follows. Introduce a set of functions 
\begin{equation*}
\mathcal{H} = \{h: \mathbb{R}_{+}\to \mathbb{R}_{+}:\, h(x)>h(y)\,\, \forall\, x>y;\,\,
\lim\limits_{x\to \infty} h(x)/x = 0\}, 
\end{equation*} 
where $\mathbb{R}_{+} = (0, +\infty)$. Then we say that $\omega$ satisfies condition $(A)$  with a function $h\in \mathcal{H}$ (see \cite{Iva21}) if there exist a subsequence 
$\{q_{n_{j}}\}_{j=1}^{\infty}$ and positive constants $C_{t}$, $C_{\delta}$ such that $C_{\delta}<1$
\begin{equation}\label{eqCondA1}
q_{n_{j}+1} >  q_{n_{j}} h(q_{n_{j}}),\;\; \forall\; j\in \mathbb{N}
\end{equation}
and for all $k\in \mathbb{N}$ there exists an index $J_{k}$:
\begin{equation}\label{eqCondA2}
\frac{1}{q_{n_{J_{k}}}}\left(\log\left( q_{n_{J_{k}}}h( q_{n_{J_{k}}})\right)+ 
\log C_{\delta}^{-1}\right)<
\frac{\log \left(q_{n_{J_{k+1}}} h(q_{n_{J_{k+1}}})\right)}{q_{n_{J_{k}}}} < C_{t}.
\end{equation}
It has to be noted that the set of those $\omega\in (0, 1)$ which satisfy the Brjuno's condition is of full measure  \cite{Brj}. The situation with the condition $(A)$ depends on a function 
$h\in \mathcal{H}$. It is known (see e.g. \cite{Khi})  that, on the one hand, for all 
$\omega\in (0, 1)$ the denominators $q_{n}$ satisfy
\begin{equation*}
q_{n} > 2^{\frac{n-1}{2}},\;\; \forall \; n\in \mathbb{N},
\end{equation*}
whereas, on the other hand,  for almost all $\omega\in (0,1)$ there exists constant $C_{L}$ such
\begin{equation*}
q_{n} < {\rm e}^{C_{L}n},\;\; \forall \; n\in \mathbb{N}.
\end{equation*}
If the series $\sum\limits_{n=1}^{\infty} 1/h({\rm e}^{C_{L}n})$ diverges, the set of those $\omega\in (0, 1)$ which satisfy condition $(A)$ has full measure \cite{Iva21}. Particularly, 
$h(x) = \log(1 + x)$ provides an example of such functions. On the other hand, if the series $\sum\limits_{n=1}^{\infty} 1/h(2^{\frac{n-1}{2}})$ converges, 
the set of those $\omega\in (0, 1)$ which satisfy condition $(A)$ has measure zero. However, for any $h\in \mathcal{H}$ this set is dense in $(0, 1)$ \cite{Iva21}.

Under these two assumptions on the rotation number, the following theorem can be proved
\newtheorem{thms}{Theorem}
\begin{thms} 
Assume hypotheses $(H_{1}) - (H_{6})$ hold true. If $\omega$ satisfies the Brjuno's condition and the condition $(A)$ with
$$C_{t} = \log \lambda_{0} - C_{B}, $$
then there exist sufficiently small $\varepsilon_{0}>0$, positive constants $C_{\Lambda}<1$, $C_{0}$  and a subset $\mathcal{E}_{h} \subset [a, b]$ such that 
\newline 
\newline
1. the Lebesgue measure 
$leb\left([a, b]\setminus \mathcal{E}_{h}\right) = 
O\left({\rm e}^{-C_{0}/\varepsilon_{0}}\right)$;
\newline 
\newline
2. for any $t\in\mathcal{E}_{h}$ the integrated Lyapunov exponent $\hat\Lambda_{0} > C_{\Lambda}\log\lambda_{0}$;
\newline 
\newline
3. for any $t\in  \mathcal{E}_{h}$ the cocycle $M$ is NUH. 
\end{thms}

{\bf Remark 2} The proof of this theorem is similar to the proof of Theorm 4 in \cite{Iva21}. Indeed, the Brjuno's condition provides the lower bound for the Lyapunov exponent on a set of positive Lebesgue measure and, hence,  the positiveness of $\tilde\Lambda_{0}$ . On the other hand, the condition $(A)$ implies the dominance of the primary collisions for a sequence 
$\{\delta_{j}\}_{j=0}^{\infty}$ with
\begin{equation*}
\delta_{j} = \frac{1}{q_{n_{j}}h(q_{n_{j}})}.
\end{equation*}
This leads to existence of a limit critical set such that in it's neighbourhood there exist points with an arbitrary small Lyapunov exponent. Such effect can be illustrated in the following way. Consider a trajectory of a point $x$ from a small neighbourhood $U_{\delta}(c_{0})$. After $\tau_{0}(\delta)$ iterations it enters this neighbourhood again and the norm of the corresponding cocycle changes its growth. If one performs the next $\tau_{0}(\delta)$ iterations, the trajectory hits $U_{\delta}(c_{0})$ once more. However, since the products of matrices, corresponding to such two passages from $U_{\delta}(c_{0})$ to itself, are very close to each other, the norm of their superposition becomes of order $O(1)$. Finally, we remark that, to make the number of such repetitions larger, one needs to take the size of neighbourhoods, $\delta$, smaller. That is why the condition $(A)$ was imposed.

It is to be noted that, due to specific parameter dependence of the cocycle in \cite{Iva21}, the number of points, constituting the critical set, grew up with decreasing of $\varepsilon$.  To overcome this difficulty Theorem 4 in \cite{Iva21} was proved under an assumption that the rotation number satisfies the condition $(A)$  with the function $h(x) = x^{\gamma}, 0 < \gamma < 1$. As it was mentioned above, the set of such $\omega\in (0,1)$ is dense, but has measure zero. Hypotheses $(H_{1} - H_{4})$ enable us to avoid this restriction and state the result in Theorem 1 of the present paper for any function $h\in \mathcal{H}$.

\section{Secondary collisions}

In this section we study the case corresponding to dominance of the secondary collisions. As it was mentioned above, the secondary collisions depend on the parameter $t$ non-trivially. It means that, on the one hand, one may detune the parameter to avoid some specific secondary collisions and, thus, to consider only generic ones. But, on the other hand, if such generic collision occurs, it exists for some interval of $t$.

To simplify the exposition, 
we consider the simplest skew product of type (\ref{eqA}), whose critical set consists of two points
\begin{equation*}
\mathcal{C}_{0} = \{c_{0}, c_{1}\}.
\end{equation*}
We fix a positive $\delta\ll 1$ and assume that for some $n\in \mathbb{N}$
\begin{equation}\label{eqRes}
\tau_{0,1}(\delta) = n,\quad \tau_{0}(\delta) > n.
\end{equation}

We formulate two statements which are direct consequences of Lemma 1.
\begin{lems} 
Let $x\in \mathbb{T}^{1}$ be such that 
\begin{equation*}
\left\vert\cos\left(\varphi\left(\sigma_{\omega}^{k}(x)\right)\right)\right\vert\ge C_{3},\; k=0,...,n-1.
\end{equation*}
Assume also that $\lambda_{0}\gg 1$. Then there exists a positive constant $C_{M} < C_{3}$ such that
\begin{equation*}
\Vert M(x,n)\Vert \ge \left(C_{M}\lambda_{0}\right)^{n}.
\end{equation*}
\end{lems}
Thus, if the finite trajectory of a point $x\in \mathbb{T}^{1}$ does not enter a small neighbourhood of the critical set, then the cocycle $M$ satisfies condition (\ref{eqUH}).
The second statement is
\begin{lems} 
Assume $\lambda_{1}, \lambda_{2}\gg 1$. Then there exists a positive constant $C$ such that for any $\phi_{1}$ the norm, $\mu$,  of the product 
$P_{1}=Z(\lambda_{2})R(\phi_{1})Z(\lambda_{1})$ admits an estimate
\begin{align*}
&1.\quad \lambda_{1}\gg \lambda_{2} \Rightarrow \mu\ge C\frac{\lambda_{1}}{\lambda_{2}}\left(1 + O\left(\frac{\lambda_{2}^{2}}{\lambda_{1}^{2}}\right)\right),\,\, 
C < \frac{3}{2\sqrt{2}},\\
&2.\quad \lambda_{2}\gg \lambda_{1} \Rightarrow \mu\ge C \frac{\lambda_{2}}{\lambda_{1}}\left(1 + O\left(\frac{\lambda_{1}^{2}}{\lambda_{2}^{2}}\right)\right),\,\, 
C < 1.
\end{align*}
\end{lems} 

PROOF:  First, we introduce a parameter $\xi = \cos 2\phi_{1} + 1$ and a function $F = a(1+c)$, where $a$ and $c$ are from Lemma 1. Represent 
$P_{1}=Z(\lambda_{2})R(\phi)Z(\lambda_{1})$ as 
\begin{equation*}
P_{1} = R(\Phi_{1})Z(\mu_{1})R(\chi_{1}). 
\end{equation*}
Then, by Lemma 1, the norm, $\mu_{1}$, of the product $P_{1}$ is of the form
\begin{equation}\label{eq_mu_F}
\mu = \frac{1}{2}\left(F + \left( F^{2} - 4\right)^{1/2}\right)
\end{equation}
and can be considered as a function of $\xi$. Moreover, since
\begin{equation*}
F = a + \frac{1}{a} + \frac{b^{2}}{a} \ge 2,
\end{equation*}
the functions $\mu_{1}$ and $F$ (as functions of $\xi$) achieve their minima simultaneousely.

Define the following polynomials
\begin{align*}
&P(\xi) = 2\beta + (\lambda_{2}^{2} - \beta)\xi,\quad
S(\xi) = 2\beta^{2} + (1 - \beta^{2})\xi,\\
&Q(\xi) = 4\left(\frac{\lambda_{2}}{\lambda_{1}}\right)^{2}\beta^{2} +
2\left [\left(\frac{\lambda_{2}}{\lambda_{1}}\right)^{4}\gamma^{2} +
\left(\frac{\lambda_{2}}{\lambda_{1}}\right)^{2}(1 - \beta^{2})\right ]\xi - 
\left(\frac{\lambda_{2}}{\lambda_{1}}\right)^{4}\gamma^{2},\\
&R(\xi) = P^{2}(\xi) + Q(\xi),\quad 
\gamma = \frac{1 - \lambda_{2}^{-4}}{1 + \lambda_{1}^{-2}\lambda_{2}^{-2}}.
\end{align*}
Then one may represent
\begin{align*}
&F = \frac{\lambda_{1}}{\sqrt{2}\lambda_{2}} \frac{R(\xi)}{P(\xi)S^{1/2}(\xi)},\\
&F' = \frac{\lambda_{1}}{\sqrt{2}\lambda_{2}} 
\frac{P(\xi)S(\xi)R'(\xi) - R(\xi)S(\xi)P'(\xi) - \frac{1}{2}P(\xi)R(\xi)S'(\xi)}
{P^{2}(\xi)S^{3/2}(\xi)},
\end{align*}
where the prime stands for the derivative with respect to $\xi$.

Note that $PSR' - RSP' -\frac{1}{2}PRS'$ is a polynomial of degree $3$. We denote coefficients of this polynomial by $A_{k}, k=0,\ldots, 3$ and represent
\begin{equation*}
PSR' - RSP' -\frac{1}{2}PRS' = A_{0} + A_{1}\xi + A_{2}\xi^{2} + A_{3}\xi^{3}.
\end{equation*}
Taking into account definition of $P, Q, R, S$ and assumption $\lambda_{1,2}\gg 1$, one obtains
\begin{align*}
&A_{0} = 4\lambda_{2}^{2}\left(\left(\lambda_{1}^{-2} + \lambda_{2}^{-2}\right)^{4} + 
O\left(\lambda_{1}^{-12} + \lambda_{2}^{-12}\right)\right),\\
&A_{1} = -2\lambda_{2}^{4}\left(2\lambda_{1}^{-6} + 5\lambda_{1}^{-4}\lambda_{2}^{-2} + 4\lambda_{1}^{-2}\lambda_{2}^{-4} + \lambda_{2}^{-6} +
O\left(\lambda_{1}^{-10} + \lambda_{2}^{-10}\right)\right),\\
&A_{2} = -\lambda_{2}^{6}\left(\lambda_{1}^{-4} + 3\lambda_{1}^{-2}\lambda_{2}^{-2} + 2\lambda_{2}^{-4} + O\left(\lambda_{1}^{-8} + \lambda_{2}^{-8}\right)\right),\\
&A_{3} = \lambda_{2}^{6}\left(\lambda_{1}^{-4} + 3\lambda_{1}^{-2}\lambda_{2}^{-2} + 2\lambda_{2}^{-4} + O\left(\lambda_{1}^{-8} + \lambda_{2}^{-8}\right)\right).
\end{align*}
Hence, the function $F$ achieves  its minimum at $\xi_{0}\approx -A_{0}/A_{1}$. Precisely, in the case $\lambda_{2}\gg \lambda_{1}$ we have
\begin{equation*}
\xi_{0} = \lambda_{1}^{-2}\lambda_{2}^{-2}
\left(1 + O\left(\left(\frac{\lambda_{1}}{\lambda_{2}}\right)^{2}\right)\right), \quad
F(\xi_{0}) = \frac{\lambda_{2}}{\lambda_{1}}
\left(1 + O\left(\left(\frac{\lambda_{1}}{\lambda_{2}}\right)^{2}\right)\right).
\end{equation*}
 In the case $\lambda_{1}\gg \lambda_{2}$ one has
\begin{equation*}
\xi_{0} = 2\lambda_{2}^{-4}
\left(1 + O\left(\left(\frac{\lambda_{2}}{\lambda_{1}}\right)^{2}\right)\right), \quad
F(\xi_{0}) = \frac{3}{2\sqrt{2}}\frac{\lambda_{1}}{\lambda_{2}}
\left(1 + O\left(\left(\frac{\lambda_{2}}{\lambda_{1}}\right)^{2}\right)\right).
\end{equation*}
Thus, in both cases $F\gg 1$. Taking this into account together with formula (\ref{eq_mu_F}), we conclude that
\begin{equation*}
\mu_{1} = F + O(F^{-1}),
\end{equation*}
what finishes the proof. $\square$

Consider a product $P_{2} = R(\phi_{2}) Z(\lambda_{2}) R(\phi_{1}) Z(\lambda_{1})$. Then, by Lemma 1, it can be represented in a form
\begin{equation}\label{eq_P2_def}
P_{2} = 
R(\Phi_{2})Z(\mu_{2}) R(\chi_{2}),
\end{equation}
where $\Phi_{2} = \phi_{2} + \Phi_{1}$, $\Phi_{1} = \psi_{1} - \chi_{1}$, $\mu_{2} = \mu_{1}$, $\chi_{2} = \chi_{1}$ and angles 
$\psi_{1}, \chi_{1}$ together with $\mu_{1}$ are described by (\ref{eq_Lem1}).

First, we remark that, due to Lemma 3, the difference in magnitude of $\lambda_{1}$ and 
$\lambda_{2}$ (under assumption $\lambda_{1,2}\gg 1$) implies the norm of $P_{2}$ becomes large for any $\phi_{1,2}$. On the other hand, if both angles $\phi_{1,2}$ are close to $\pi/2 \mod \pi$, the sum $\phi_{2}+\Phi_{1}$ is close to $0\mod \pi$. 

Our goal now is to investigate how interaction between two small $\delta$-neighborhoods of the critical points $c_{0}, c_{1}$ influences the parameters of the cocycle $M$. To model such interaction we note that, taking $\varepsilon, \delta$ to be sufficiently small, function 
$\varphi$ can be approximated in the $\delta$-neighborhood of $c_{j}$ by
\begin{align}\nonumber
&\hat\varphi(y; L_{-}, L_{+}, k) = \pm\frac{\pi}{2} + \vartheta(y; L_{-}, L_{+}, k),\\
\label{eq_StepFunc}
&\vartheta(y; L_{-}, L_{+}, k) = 
ky\Theta(ky-L_{-})\Theta(L_{+}-ky) + L_{-}\Theta(L_{-}-ky) + L_{+}\Theta(ky-L_{+}),
\end{align}
where $y=x-c_{j}$, $\Theta$ is the Heaviside function and $k = \varphi'(c_{j})$. Note that hypotheses $(H_{2})$, $(H_{4})$ imply $L_{-}\cdot L_{+}<0$, $k \sim \varepsilon^{-1}$. A graph of such function $\hat\varphi$ is shawn on Fig.1. 

\begin{figure}[h]
\center{\includegraphics[width=0.4\linewidth]{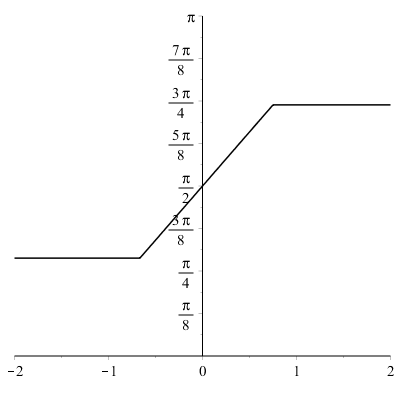}}
\caption{Graph of function $\hat\varphi$ with $L_{-} = -2/3$, $L_{+} = 3/4$ and $k=1$.}
\end{figure}

More precisely, hypotheses ($H_{2}$), ($H_{3}$) guarantee the existence of a positive $\delta > C_{3}\varepsilon$ such that 
\begin{equation}\label{eq_phi_approx}
\left\vert \hat\varphi(x-c_{j};  L_{-,\min}^{(j)}, L_{+,\min}^{(j)}, k^{(j)})\right\vert \le 
\left\vert \varphi(x)\right\vert \le 
\left\vert \hat\varphi(x-c_{j};  L_{-,\max}^{(j)}, L_{+,\max}^{(j)}, k^{(j)})\right\vert,\quad \forall\, x\in U_{\delta}(c_{j})
\end{equation}
with some parameters $L_{\pm,\min}^{(j)}, L_{\pm,\max}^{(j)}, k^{(j)}$.
Taking this into account we set the angles $\phi_{j}$ in the product $P_{2}$ to be equal to
\begin{equation}\label{eq_Inter}
\phi_{1}(x) = \hat\varphi(x-\Delta; L_{-}^{(1)}, L_{+}^{(1)}, k^{(1)}),\quad
\phi_{2}(x) = \hat\varphi(x; L_{-}^{(2)}, L_{+}^{(2)}, k^{(2)})
\end{equation} 
with fixed parameters $L_{\pm}^{(j)}, k^{(j)} = \varepsilon^{-1}r^{(j)}, j=1,2$ and a detuning parameter 
$\Delta$, which will be specified later. Due to  assumption $(H_{4})$, one has 
 $\varphi(c_{0}) = \varphi(c_{1}) \mod 2\pi$. Hence, the signs before $\pi/2$ in the definition of $\phi_{1,2}$ coincide. We consider the case of sign '+' (the case of sign '-' can be studied similarly). 

\begin{lems} 
Assume $\lambda_{1}, \lambda_{2}\gg 1$ and $\lambda_{1}\gg \lambda_{2}$. 
If $\phi_{1}, \phi_{2}$ are described by (\ref{eq_Inter}) with $\Delta$ satisfying
\begin{equation*}
\vert\Delta\vert < \frac{\varepsilon \sqrt{\beta}}{\left(r^{(1)} r^{(2)}\right)^{1/2}},
\end{equation*}
then product $P_{2}$ has the following characteristics:
\begin{align*}
&\mu_{2} \ge \frac{\lambda_{1}}{\lambda_{2}}
\left(1 + O\left(\frac{\lambda_{2}^{2}}{\lambda_{1}^{2}}\right)\right),\\
&\vert\chi_{2}\vert_{\max} = 
\frac{\gamma}{2}\left(\frac{\lambda_{2}}{\lambda_{1}}\right)^{2}
\left(1 + O\left(\frac{1}{\lambda_{2}^{4}} + 
\frac{\lambda_{2}^{2}}{\lambda_{1}^{2}}\right)\right),\\
&\vert\Phi_{2}\vert_{\max} = \frac{\pi}{2} - 
\sqrt{\beta}\left\vert\frac{r^{(2)}}{r^{(1)}}\right\vert^{1/2}
\left(1 + O\left(\frac{1}{\lambda_{2}^{4}} + 
\frac{\lambda_{2}^{2}}{\lambda_{1}^{2}}\right)\right),
\end{align*}
where $\vert\chi_{2}\vert_{\max}, \vert\Phi_{2}\vert_{\max}$ stand for the maximum values of $\chi_{2}$ and $\Phi_{2}$, respectively.
\end{lems} 

PROOF: Substituting (\ref{eq_Inter}) into (\ref{eq_Lem1}) one obtains expressions for $\psi_{1}, \chi_{1}$ as functions of variable $x$. It is to be noted that $\vert \tan\chi\vert$ is an increasing, odd function of a variable $\eta = \sqrt{2}b(1-c)^{-1}$. On the other hand, using notations from Lemma 3, we have
\begin{equation}\label{eq_eta}
\vert\eta\vert = \varkappa\gamma\frac{\sqrt{\xi(2-\xi)}P(\xi)}{P^{2}(\xi)-Q(\xi)},\quad    \varkappa = \left(\frac{\lambda_{2}}{\lambda_{1}}\right)^{2}.
\end{equation} 
Differentiating (\ref{eq_eta}) with respect to $\xi$, one concludes that $\eta$ achieves maximum at $\xi_{*}$, which solves the following equation
\begin{equation}\label{eq_xi_star}
\xi(2-\xi)\left(P(\xi)Q'(\xi) - \left(P^{2}(\xi) + Q'(\xi)\right)P'(\xi)\right) + 
(1-\xi)P(\xi)\left(P^{2}(\xi) - Q(\xi)\right) = 0.
\end{equation}
We note that the left hand side of (\ref{eq_xi_star}) is a polynomial of degree 3 and can be represented as
\begin{equation*}
\xi(2-\xi)\left(P Q' - \left(P^{2} + Q'\right)P'\right) + 
(1-\xi)P\left(P^{2} - Q\right) = B_{0} + B_{1}\xi + B_{2}\xi^{2} + B_{3}\xi^{3},
\end{equation*}
where $B_{k}, k=0,\ldots, 3$ are constants. Under assumptions $\lambda_{1,2}\gg 1$, $\varkappa\ll 1$ the coefficients $B_{j}$ admit the following representation
\begin{align*}
&B_{0} = 8\beta^{3}\left(1 + O\left(\beta \lambda_{2}^{-2} + \varkappa\right)\right),\\
&B_{1} = 4\beta^{2}\lambda_{2}^{2}
\left(1 + O\left(\beta \lambda_{2}^{-2} + \varkappa\right)\right),\\
&B_{2} = -2\beta\lambda_{2}^{4}
\left(1 + O\left(\beta \lambda_{2}^{-2} + \varkappa\right)\right),\\
&B_{3} = -\lambda_{2}^{6}
\left(1 + O\left(\beta \lambda_{2}^{-2} + \varkappa\right)\right).
\end{align*}
The only positive solution of (\ref{eq_xi_star}) satisfies
\begin{equation*}
\xi_{*} = 2\beta\lambda_{2}^{-2}
\left(1 + O\left(\beta \lambda_{2}^{-2} + \varkappa\right)\right).
\end{equation*}
and the maximum value of $\eta$ equals
\begin{equation}\label{eq_eta_max}
\eta_{\max} = \eta(\xi_{*}) = \frac{\varkappa \gamma}{2\lambda_{2}\sqrt{\beta}}
\left(1 + O\left(\beta \lambda_{2}^{-2} + \varkappa\right)\right).
\end{equation}
Thus, as $\lambda_{1,2}\gg 1$ and $\varkappa\ll 1$ the following estimate holds
\begin{equation}\label{eq_tan_chi}
\tan \chi = \sqrt{2} b \left(1 + O\left(\varkappa\right)\right) = 
\varkappa\gamma\frac{\sin\phi \cos\phi}
{\lambda_{2}^{2}\cos^{2}\phi + \beta\sin^{2}\phi}
\left(1 + O\left(\varkappa\right)\right).
\end{equation}
The graph of $\tan \chi$ as a function of $\phi$ is presented on Fig. 2.

\begin{figure}[h]
\center{\includegraphics[width=0.4\linewidth]{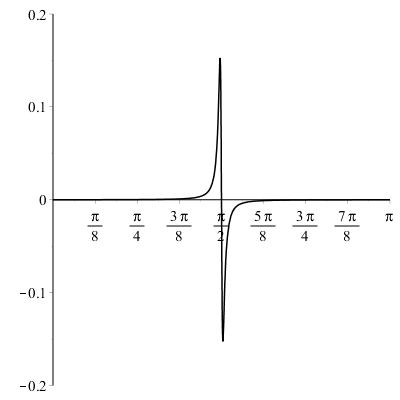}}
\caption{Graph of $\tan \chi$, corresponding to $\lambda_{1} = 10, \lambda_{2}=7$.}
\end{figure}

We substitute (\ref{eq_Inter}) into (\ref{eq_Lem1}) and conclude that under assumptions
$\lambda_{1,2}\gg 1$, $\varkappa\ll 1$, $\varepsilon\ll 1$
\begin{equation}\label{eq_Phi1}
\tan \Phi_{1}(x) = \alpha_{1}(x)\tan \phi_{1}(x),\quad 
\alpha_{1}(x) = \frac{\beta  + 
\varkappa\gamma\frac{\cos^{2}\phi_{1}(x)}
{\lambda_{2}^{2}\cos^{2}\phi_{1}(x) + \beta\sin^{2}\phi_{1}(x)}
\left(1 + O\left(\varkappa\right)\right)}
{1 - \beta\varkappa\gamma\frac{\sin^{2}\phi_{1}(x)}
{\lambda_{2}^{2}\cos^{2}\phi_{1}(x) + \beta\sin^{2}\phi_{1}(x)}
\left(1 + O\left(\varkappa\right)\right)},
\end{equation}
where the factor $\alpha_{1}$ satisfies 
\begin{equation}\label{eq_alpha_1}
\beta < \alpha_{1}(x) < 
\frac{\beta + \varkappa\gamma\lambda_{2}^{-2}}{1-\varkappa\gamma} = 
\beta \left(1 + O\left(\varkappa\right)\right).
\end{equation}
Hence, the graph of $\Phi_{1}(x)$ is very similar to those one of $\phi_{1}(x)$ exept the slope of the graph at $x=\Delta$ is much higher for $\Phi_{1}(x)$  than for $\phi_{1}(x)$, since 
$\Phi'_{1}(0) = \varepsilon^{-1}\beta^{-1}r^{(1)}\left(1 + O\left(\varkappa\right)\right)$. Moreover, new levels $\hat L_{\pm} = \arctan \left(\alpha_{1}(\pm\delta)\tan\left(L_{\pm}\right)\right)$ become close to $0 \mod \pi$ as $\alpha\ll 1$. Particularly, one has
\begin{equation}\label{eq_L}
\hat L_{\pm} = O(\beta).
\end{equation}

We represent $\Phi_{1}(x)$ as
\begin{equation*}
\Phi_{1}(x) = \frac{\pi}{2} + 
\arctan\left(\alpha_{1}^{-1}(x)\tan
\left(\vartheta(x-\Delta, L_{-}^{(1)}, L_{+}^{(1)}, k^{(1)})\right)\right) 
\end{equation*}
and consider the sum $\Phi_{2}(x) = \phi_{2}(x) + \Phi_{1}(x)$ 
\begin{align}\nonumber
&\Phi_{2}(x) = \vartheta_{2}(x) + \hat\vartheta_{1}(x-\Delta) \mod \pi,\\
\label{eq_Psi_2}
&\vartheta_{2}(x) = \vartheta(x, L_{-}^{(2)}, L_{+}^{(2)}, k^{(2)}),\\ 
\nonumber
&\hat\vartheta_{1}(x) = \arctan\left(\alpha_{1}^{-1}(x)\tan
\left(\vartheta(x, L_{-}^{(1)}, L_{+}^{(1)}, k^{(1)})\right)\right).
\end{align}
Note that graphs of both summands in the right hand side of (\ref{eq_Psi_2}) are of the same shape. However, as we study the secondary collisions, assumption ($H_{4}$) implies 
\begin{equation*}
\varphi'(c_{0})\cdot \varphi'(c_{1}) < 0
\end{equation*}
or equivalently $k^{(1)}\cdot k^{(2)} < 0$. Fig. 3 illustrates relative position of these two graphs.

\begin{figure}[h]
\begin{minipage}[h]{0.49\linewidth}
\center{\includegraphics[width=1.2\linewidth]{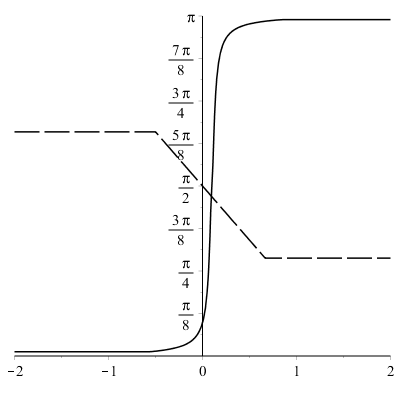} \\ a)}
\end{minipage}
\begin{minipage}[h]{0.49\linewidth}
\center{\includegraphics[width=1.2\linewidth]{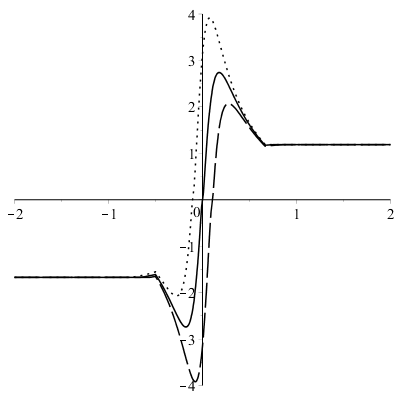} \\ b)}
\end{minipage}
\caption{a. Graphs of functions $\Phi_{1}$ (solid line) and $\phi_{2}$ (dashed line), corresponding to $\lambda_{1}=10, \lambda_{2} =7, L_{-}^{(1)} = -2/3, L_{+}^{(1)}=3/4, k^{(1)} = 1, 
L_{-}^{(2)} = 2/3, L_{+}^{(2)}=-1/2, k^{(2)} = -1, \Delta = 0.1$;  b. graph of $\tan \Psi_{2}$, corresponding to different values of $\Delta$ ($\Delta=0$ - solid line, $\Delta = 0.1$ - dashed line, $\Delta=-0.1$ - dotted line; values of other parameters are the same as at Fig.3a)}
\end{figure}

 If $\Delta = 0$ both $\vartheta_{2}(x)$ and $\hat\vartheta_{1}(x)$ vanish at $x=0$, but are of different signs as $x\neq 0$. Hence, differentiating (\ref{eq_Psi_2}) one obtains that $\vert\Phi_{2}\vert$ attains maximum at $x_{*}$, which satisfies
\begin{equation}\label{eq_x_star}
\alpha^{-1}(x)\vert k^{(1)}\vert - \vert k^{(2)}\vert = 
\left(\alpha^{-2}(x) - 1\right) \vert k^{(2)}\vert \sin^{2}\left(k^{(1)} x\right).
\end{equation}
Assumption ($H_{2}$) together with (\ref{eq_alpha_1}) yields 
\begin{equation*}
x_{*} = \frac{1}{k^{(1)}} \left\vert\frac{k^{(1)}}{k^{(2)}}\right\vert^{1/2} 
\left(1 + O(\beta)\right)
\end{equation*}
and
\begin{equation*}
\vert\Phi_{2}\vert_{\max} = \vert\Phi_{2}(x_{*})\vert = \frac{\pi}{2} - 2\sqrt{\beta}\left\vert\frac{k^{(2)}}{k^{(1)}}\right\vert^{1/2}\left(1 + O(\beta)\right).
\end{equation*}
If $\Delta\neq 0$ then $\vert\Phi_{2}\vert$ attains maximum at $x_{*}$
\begin{equation*}
x_{*} = \Delta - \textrm{sign}(\Delta)\frac{1}{k^{(1)}} \left\vert\frac{k^{(1)}}{k^{(2)}}\right\vert^{1/2} 
\left(1 + O(\beta)\right)
\end{equation*}
and
\begin{equation*}
\vert\Phi_{2}\vert_{\max} = \vert\Phi_{2}(x_{*})\vert = \frac{\pi}{2} - 2\sqrt{\beta}\left\vert\frac{k^{(2)}}{k^{(1)}}\right\vert^{1/2}\left(1 + O(\beta)\right) + 
\frac{\Delta}{k^{(2)}}.
\end{equation*}
We assume that
\begin{equation}\label{eq_Delta}
\vert\Delta\vert < \frac{\varepsilon \sqrt{\beta}}{\left(r^{(1)} r^{(2)}\right)^{1/2}}.
\end{equation}
In this case the maximum value $\vert\Phi_{2}\vert_{\max}$ is bounded by
\begin{equation*}
\vert\Phi_{2}\vert_{\max} < \frac{\pi}{2} - \sqrt{\beta}\left\vert\frac{k^{(2)}}{k^{(1)}}\right\vert^{1/2}\left(1 + O(\beta)\right).
\end{equation*}
This finishes the proof. $\square$ 

{\bf Remark 3}
One may note that, if functions $\phi_{1,2}$ in the product $P_{2}$  satisfy 
\begin{align}\label{eq_phi12}
\nonumber
&\left\vert \hat\varphi(x-\Delta;L_{-,\min}^{(1)}, L_{+,\min}^{(1)},k_{\min}^{(1)})\right\vert \le 
\left\vert \phi_{1}(x)\right\vert \le 
\left\vert \hat\varphi(x-\Delta;L_{-,\max}^{(1)}, L_{+,\max}^{(1)},k_{\max}^{(1)})\right\vert,\\
&\left\vert \hat\varphi(x;L_{-,\min}^{(2)}, L_{+,\min}^{(2)},k_{\min}^{(2)})\right\vert \le 
\left\vert \phi_{2}(x)\right\vert \le 
\left\vert \hat\varphi(x;L_{-,\max}^{(2)}, L_{+,\max}^{(2)},k_{\max}^{(2)})\right\vert
\end{align}
and indices $j_{1}, j_{2}\in \{1,2\}$ are such that 
\begin{equation*}
j_{1}\neq j_{2},\quad k_{\max}^{(j_{1})} \ge k_{\max}^{(j_{2})},
\end{equation*}
then
\begin{align*}
&\left\vert \phi_{1}(x) - \phi_{2}(x)\right\vert \le 
\left\vert \vartheta(x-\Delta;L_{-,\max}^{(1)}, L_{+,\max}^{(1)},k_{\max}^{(1)}) -
\vartheta(x;L_{-,\min}^{(2)}, L_{+,\min}^{(2)},k_{\min}^{(2)})\right\vert,\quad j_{1} = 1;\\
&\left\vert \phi_{1}(x) - \phi_{2}(x)\right\vert \le 
\left\vert \vartheta(x-\Delta;L_{-,\min}^{(1)}, L_{+,\min}^{(1)},k_{\min}^{(1)}) -
\vartheta(x;L_{-,\max}^{(2)}, L_{+,\max}^{(2)},k_{\max}^{(2)})\right\vert,\quad j_{1} = 2.
\end{align*}
Taking this into account we arrive at the following
\newtheorem{cors}{Corollary}
\begin{cors}
Assume $\lambda_{1}, \lambda_{2}\gg 1$, $\lambda_{1}\gg \lambda_{2}$ and 
$\phi_{1}, \phi_{2}$ satisfy (\ref{eq_phi12}) with $\Delta$ such that
\begin{equation*}
\vert\Delta\vert < \frac{\varepsilon \sqrt{\beta}}{\left(r_{\max}^{(j_{1})} r_{\min}^{(j_{2})}\right)^{1/2}}.
\end{equation*}
Then product $P_{2}$ has the following characteristics:
\begin{align*}
&\mu_{1} \ge \frac{\lambda_{1}}{\lambda_{2}}
\left(1 + O\left(\frac{\lambda_{2}^{2}}{\lambda_{1}^{2}}\right)\right),\\
&\vert\chi_{1}\vert_{\max} \le 
\frac{\gamma}{2}\left(\frac{\lambda_{2}}{\lambda_{1}}\right)^{2}
\left(1 + O\left(\frac{1}{\lambda_{2}^{4}} + 
\frac{\lambda_{2}^{2}}{\lambda_{1}^{2}}\right)\right),\\
&\vert\Phi_{2}\vert_{\max} \le \frac{\pi}{2} - 
\sqrt{\beta}\left\vert\frac{r_{\min}^{(j_{2})}}{r_{\max}^{(j_{1})}}\right\vert^{1/2}
\left(1 + O\left(\frac{1}{\lambda_{2}^{4}} + 
\frac{\lambda_{2}^{2}}{\lambda_{1}^{2}}\right)\right).
\end{align*}
\end{cors}

Finally, consider a product 
$P_{3} = Z(\lambda_{3})P_{2} = 
Z(\lambda_{3})R(\phi_{2})Z(\lambda_{2})R(\phi_{1})Z(\lambda_{1})$. By Lemma 1 we may represent it in a form
\begin{equation}\label{eq_P3}
P_{3} = R(\Phi_{3})Z(\mu_{3})R(\chi_{3}).
\end{equation}
Then we obtain
\begin{lems} 
Assume conditions of Lemma 4 are satisfied. If, additionally, 
$\lambda_{1}\gg \lambda_{2}^{3/2}$ and $\lambda_{3}>\lambda_{2}^{1/2}$, 
then product $P_{3}$ has the following characteristics:
\begin{align*}
&\left\vert\tan\Phi_{3}\right\vert <  
\left\vert\frac{r^{(1)}}{r^{(2)}}\right\vert^{1/2}
\left(1 + O\left(\frac{1}{\lambda_{2}^{4}} + 
\frac{\lambda_{2}^{2}}{\lambda_{1}^{2}}\right)\right),\\
&\mu_{3} \ge \frac{1}{2}\left(1 + \frac{r^{(2)}}{r^{(1)}}\right)
\frac{\lambda_{1}}{\lambda_{2}^{3/2}}
\left(1 + O\left(\frac{1}{\lambda_{2}^{4}} + 
\frac{\lambda_{2}^{3}}{\lambda_{1}^{2}}\right)\right).
\end{align*}
\end{lems} 

Proof: Using representation (\ref{eq_P2_def}), we apply Lemmas 1 and 4 to the product
$Z(\lambda_{3})R(\Phi_{2})Z(\mu_{2})$ and obtain that
\begin{equation*}
Z(\lambda_{3})R(\Phi_{2})Z(\mu_{2}) = R(\Phi_{3})Z(\mu_{3})R(\chi_{3}-\chi_{2})
\end{equation*}
with
\begin{equation*}
\left\vert \tan\Phi_{3}\right\vert < 
\beta_{2}\sqrt{\beta}\left\vert\frac{r^{(1)}}{r^{(2)}}\right\vert^{1/2}
\left(1 + O\left(\beta_{2}\lambda_{3}^{-2} + \beta\lambda_{2}^{-2} + \varkappa\right)\right),\quad
\beta_{2} = \frac{\lambda_{3}^{-2} + \mu_{1}^{-2}}
{1 + \lambda_{3}^{-2} \mu_{1}^{-2}}. 
\end{equation*}
Then, taking into account that $\beta = \lambda_{2}^{-2}(1+O(\varkappa))$ as 
$\varkappa\ll 1$ and assuming $\lambda_{3}>\lambda_{2}^{1/2}$, one obtains
\begin{equation*}
\left\vert \tan\Phi_{3}\right\vert < 
\left\vert\frac{r^{(1)}}{r^{(2)}}\right\vert^{1/2}
\left(1 + O\left(\lambda_{2}^{-4} + \varkappa\right)\right). 
\end{equation*}
Besides, Lemmas 1 and 4 imply that parameter $\mu_{3}$ (see (\ref{eq_Lem1}) for definition), corresponding to product $P_{3}$, admits an estimate
\begin{equation}\label{eq_a_min}
\mu_{3} > \frac{\mu_{1}}{\lambda_{3}}\cdot 
\frac{\left(\frac{\lambda_{3}}{\lambda_{2}}\right)^{2}\frac{r^{(2)}}{r^{(1)}} + \lambda_{3}^{-2} + \mu_{1}^{-2}}
{\left(\lambda_{2}^{-2}\frac{r^{(2)}}{r^{(1)}} + 
\left(\lambda_{3}^{-2} + \mu_{1}^{-2}\right)^{2}\right)^{1/2}}
\left(1 + O\left(\beta + \varkappa\right)\right).
\end{equation}
An expression in the right hand side of (\ref{eq_a_min}) attains its minimum with respect to $\lambda_{3}$ at $\lambda_{3,*}$, which is a solution of
\begin{equation}
3\frac{r^{(2)}}{r^{(1)}}\left(\mu_{1}^{2}\lambda_{2}^{-1}\right)^{2}\lambda_{3}^{-1}+ \left(\frac{r^{(2)}}{r^{(1)}}\right)^{2}\left(\mu_{1}^{2}\lambda_{2}^{-1}\right)^{4} + \frac{r^{(2)}}{r^{(1)}}\left(\mu_{1}^{2}\lambda_{2}^{-1}\right)^{2} =
\lambda_{3}^{-1}\left(\lambda_{3}^{-1} + 1\right)^{3}.
\end{equation}
Lemma 4 and the condition $\lambda_{1} \gg \lambda_{2}^{3/2}$ yield 
$\mu_{1}^{2}\lambda_{2}^{-1}\gg 1$ and
\begin{equation*}
\lambda_{3,*} = \left(\frac{r^{(1)}}{r^{(2)}}\right)\mu_{1}^{2}\lambda_{2}^{-1}
\left(1 + O\left(\beta + \varkappa + \mu_{1}^{-2}\lambda_{2}\right)\right).
\end{equation*}
Moreover, we obtain
\begin{equation*}
\mu_{3} > \mu_{1}\lambda_{2}^{-1/2}\frac{1}{2}\left(1+ \frac{r^{(2)}}{r^{(1)}}\right)
\left(1 + O\left(\beta + \varkappa + \mu_{1}^{-2}\lambda_{2}\right)\right) \gg 1.
\end{equation*}
This finishes the proof. $\square$

We apply Lemmas 1 - 5 to the cocycle $M$. Assuming that for a positive $\delta\ll 1$ the secondary collision occures, i.e.
\begin{equation*}
\tau_{0,1}(\delta) = n,\quad \tau_{0}(\delta) > n,\quad n\in \mathbb{N}, 
\end{equation*}
one may note that finite trajectory $\sigma_{\omega}^{k}(x), k = 1,\ldots, n-1$ of any point $x\in U_{\delta}(c_{0})$ does not fall into $U_{\delta}(\mathcal{C}_{0})$. Then,
applying consecutively Lemma 1, we represent the product $Z(\lambda(\sigma_{\omega}^{n}(x)))M(\sigma_{\omega}(x),n-1)$ 
for $x\in U_{\delta}(c_{0})$ in the following way
\begin{equation}\label{eqZM}
Z(\lambda(\sigma_{\omega}^{n}(x)))M(\sigma_{\omega}(x),n-1) =
R(\Phi_{n}(x)) Z(\mu_{n}(x)) R(\chi_{n}(x)).
\end{equation}
Taking into account conditions (\ref{eqRes}) and ($H_{5}$), one obtains that the angles 
$\Phi_{n}(x)$, $\chi_{n}(x)$ are small, particularly, 
\begin{equation}\label{eq_Phi_chi}
\Phi_{n}(x) = O(\lambda_{0}^{-2}),\quad \chi_{n}(x) = O(\lambda_{0}^{-2}).
\end{equation}
Moreover, $\mu_{n}(x)$ is large and satisfies
\begin{equation}\label{eq_res_mu}
\left(C_{M}\lambda_{\min}\right)^{n} \le \mu_{n}(x) \le \lambda_{\max}^{n},
\end{equation}
where $\lambda_{\min}, \lambda_{\max}$ stand for the minimum and maximum values of function $\lambda$ over $x\in \mathbb{T}^{1}$ and $C_{M}$ is a constant from Lemma 2.
Besides, we suppose that there exist integers $n_{-}, n_{+}$ such that
\begin{equation}\label{eq_cond_res}
n + n_{-} + n_{+} \le \tau_{0}(\delta),\quad
\left(C_{M}\lambda_{\min}\right)^{n_{-}} > \lambda_{\max}^{3n/2},\quad
\left(C_{M}\lambda_{\min}\right)^{n_{+}} > \lambda_{\max}^{n/2}.
\end{equation}
Then we consider two products $Z(\lambda(x))M(\sigma_{\omega}^{-n_{-}}(x), n_{-})$ and 
$Z(\lambda(\sigma_{\omega}^{n+n_{+}}(x)))M(\sigma_{\omega}^{n+1}(x), n_{+}-1)$
and represent them in a form
\begin{align}\nonumber
&Z(\lambda(x))M(\sigma_{\omega}^{-n_{-}}(x), n_{-}) = 
R(\Phi_{n,-}(x)) Z(\mu_{n,-}(x)) R(\chi_{n,-}(x)),\\
\label{eq_2prod}
&Z(\lambda(\sigma_{\omega}^{n+n_{+}}(x)))M(\sigma_{\omega}^{n+1}(x), n_{+}-1) = 
R(\Phi_{n,+}(x)) Z(\mu_{n,+}(x)) R(\chi_{n,+}(x)).
\end{align}
We note that $\Phi_{n,\pm}, \chi_{n,\pm}$ satisfy (\ref{eq_Phi_chi}) and
\begin{equation*}
\left(C_{M}\lambda_{\min}\right)^{n_{\pm}} \le \mu_{n_{\pm}}(x) \le \lambda_{\max}^{n_{\pm}}.
\end{equation*}
Hence, the cocycle $M(\sigma_{\omega}^{-n_{-}}(x), n_{-} + n + n_{+})$ reads
\begin{multline}\label{eq_Res}
M(\sigma_{\omega}^{-n_{-}}(x), n_{-} + n + n_{+}) = 
R\Bigl(\Phi_{n,+}(x)\Bigr) 
Z\Bigl(\mu_{n,+}(x)\Bigr) 
R\Bigl(\chi_{n,+}(x) + \varphi(\sigma_{\omega}^{n}(x)) + \Phi_{n}(x)\Bigr) 
Z\Bigl(\mu_{n}(x)\Bigr)\times\\ 
R\Bigl(\chi_{n}(x) + \varphi(x) + \Phi_{n,-}(x)\Bigr) 
Z\Bigl(\mu_{n,-}(x)\Bigr) 
R\Bigl(\chi_{n,-}(x)\Bigr)
\end{multline}
The central part of this product, namely, 
$R\Bigl(-\Phi_{n,+}(x)\Bigr)
M\Bigl(\sigma_{\omega}^{-n_{-}}(x), n_{-} + n + n_{+}\Bigr)
R\Bigl(-\chi_{n,-}(x)\Bigr)$
has the structure of product $P_{3}$ from Lemma 5 with parameters
\begin{align}\nonumber
&\phi_{1} = \chi_{n}(x) + \varphi(x) + \Phi_{n,-}(x),\quad
\phi_{2} = \chi_{n,+}(x) + \varphi(\sigma_{\omega}^{n}(x)) + \Phi_{n}(x),\\
&\lambda_{1} = \mu_{n,-}(x),\quad
\lambda_{2} = \mu_{n}(x),\quad
\lambda_{3} = \mu_{n,+}(x).
\end{align}
We note that provided $\delta$ to be sufficiently small, $\phi_{1}(x)\approx \varphi(x)$, 
$\phi_{2}(x)\approx \varphi(\sigma_{\omega}^{n}(x))$ as $x\in U_{\delta}(c_{0})$. Particularly, there exist unique solutions $x_{1}, x_{2}\in U_{\delta}(c_{0})$ of the following equations
\begin{equation}\label{eq_x1x2}
\cos\phi_{j}(x_{j}) = 0,\quad j=1,2.
\end{equation}
These solutions admit estimates
\begin{align}\label{eq_x1x2_est}
&x_{1} = c_{0} - \frac{\varepsilon}{r_{0}}\Bigl(\chi_{n}(c_{0}) + \Phi_{n,-}(c_{0})\Bigr)
\left(1 + O(\delta)\right),\\
&x_{2} = \sigma_{\omega}^{-n}(c_{1}) - 
\frac{\varepsilon}{r_{1}}\Bigl(\chi_{n,+}\left(\sigma_{\omega}^{-n}(c_{1})\right) + 
\Phi_{n}\left(\sigma_{\omega}^{-n}(c_{1})\right)\Bigr)
\left(1 + O(\delta)\right).
\end{align}
We emphasize here that assumption ($H_{6}$) guarantees non-trivial dependence of the distance between points $c_{0}, c_{1}$ on the parameter $t$. Without loss of generality we may assume that position of $c_{0}$ is fixed and $c_{1}$ varies with respect to $t$. Then solution $x_{2}$ can be considered as a function of $t$, whereas $x_{1}$ is constant.  

\begin{defs}
We say that a value $t_{res}$ is resonant of order $n$ if it solves
\begin{equation*}
x_{2}(t) = x_{1}.
\end{equation*}
\end{defs}
It means that for resonant $t_{res}$ the $n$-th iteration $\sigma_{\omega}^{n}(c_{0})$ falls not exactly on $c_{1}$, but close to it, i.e.
\begin{equation*}
\sigma_{\omega}^{n}(c_{0}) = c_{1}(t_{res}) + \Delta_{n}^{res}(t_{res}),
\end{equation*}
where $\Delta_{n}^{res}(t)$ is a small correction such that
\begin{equation*}
\Delta_{n}^{res}(t) = 
\varepsilon\left[
r_{0}^{-1}\left(\chi_{n}(c_{0}) + \Phi_{n,-}(c_{0})\right) -
r_{1}^{-1}\left(\chi_{n,+}(\sigma_{\omega}^{-n}(c_{1}(t))) + 
\Phi_{n}(\sigma_{\omega}^{-n}(c_{1}(t)))\right)\right]
\left(1 + O(\delta)\right).
\end{equation*}
As parameters $\Phi_{n}, \Phi_{n,\pm}, \chi_{n}, \chi_{n,\pm}$ are of the order of $\lambda_{0}^{-2}$, we obtain
\begin{equation*}
\Delta_{n}^{res} = O(\varepsilon \lambda_{0}^{-2}).
\end{equation*}

{\bf Remark 4} In the simplest case $n=1$ one has
\begin{equation*}
\Delta_{n}^{res}\equiv 0.
\end{equation*}
We note that if there exists $t_{0}\in [a,b]$ such that $\sigma_{\omega}^{-n}(c_{1}(t_{0}))=c_{0}$ one may apply the implicit function theorem to conclude the existence of $t_{res}$, which satisfies
\begin{equation*}
t_{res} = t_{0} + O(\varepsilon).
\end{equation*}

Let parameter $t$ be close to a resonant one, i.e. we suppose that $\tau_{0,1}(\delta) = n$ and
\begin{equation*}
dist\left(\sigma_{\omega}^{n}(c_{0}), c_{1}\right) = \Delta_{n} \ll 1.
\end{equation*}
Consider cocycle $M(\sigma_{\omega}^{-n_{-}}(x), n_{-} + n + n_{+})$ for $x\in I_{0}(\varepsilon)$ and its singular decomposition (\ref{eq_Res}). Then we arrive at the following
\begin{lems} 
Let $n\in \mathbb{N}$ be fixed. Assume that  $\lambda_{0}\gg 1$ and there exist integers $n_{-}, n_{+}$ such that (\ref{eq_cond_res}) holds. If, additionally, there exists $t_{0}\in [a,b]$ such that  $\sigma_{\omega}^{-n}(c_{1}(t_{0}))=c_{0}$, then there exists $\varepsilon_{0}$ such that for any $\varepsilon\in (0,\varepsilon_{0})$ the following holds:
there exists unique resonant value $t_{res}$ in a small neighborhood of $t_{0}$ and
\begin{equation*}
\Delta_{n}^{res} = -\frac{\Phi_{n}(c_{0})}{\varphi'(c_{1})}\biggl(1+O(\lambda_{0}^{-2})\biggr).
\end{equation*}
Moreover, for any $x\in U_{\delta}(c_{0})$ and any $t\in (t_{res}-h, t_{res}+h)$ the cocycle $M(\sigma_{\omega}^{-n_{-}}(x), n_{-} + n + n_{+})$ admits the following representation
\begin{equation*}
M(\sigma_{\omega}^{-n_{-}}(x), n_{-} + n + n_{+}) = \prod\limits_{k=-n_{-}}^{n+n_{+}}A(\sigma_{\omega}^{k}(x)) =\\
 R(\psi_{n,n_{-},n_{+}}(x)) Z(\mu_{n,n_{-},n_{+}}(x)) R(\chi_{n,n_{-},n_{+}}(x))
\end{equation*}
such that
\begin{equation*}
\mu_{n,n_{-},n_{+}}(x) \ge \lambda_{0}, \quad
\psi_{n,n_{-},n_{+}}(x) = \varphi(\sigma_{\omega}^{n+n_{+}}(x)) + O(\lambda_{0}^{-2}),\quad
\chi_{n,n_{-},n_{+}}(x) = O(\lambda_{0}^{-2})
\end{equation*}
and the size of a neighborhood of $t_{res}$ possesses an estimate
\begin{equation}\label{eq_lm6_h}
h = O\Bigl(\varepsilon \lambda_{0}^{-2n} \Bigr).
\end{equation}
\end{lems}

PROOF:  The proof follows from the direct application of Lemmas 1-6 to the presentation (\ref{eq_Res}) of the cocycle $M(\sigma_{\omega}^{-n_{-}}(x), n_{-} + n + n_{+})$. We only mention that to due to smooth dependence of all objects on the parameter $t$ and Lemma 4, the size of a neighborhood, $h$, can be bounded by
\begin{equation*}
h = O\Bigl(\varepsilon^{2}\mu_{n}^{-2}\bigl(c_{0}(t_{res})\bigr)\Bigr),
\end{equation*}
where $\mu_{n}(x)$ is defined by (\ref{eqZM}). Then estimate (\ref{eq_lm6_h}) is a consequence of (\ref{eq_res_mu}). $\square$

Lemma 6 is a key tool for establishing the hyperbolicity of the cocycle $M$. First, we introduce the following notations. We consider the matrix $A(x)$ defined by (\ref{eqA}) and denote it by letter $B$ if $x\in U_{\delta}(\mathcal{C}_{0})$ . In the case $x\notin U_{\delta}(\mathcal{C}_{0})$, we denote $A(x)$ by letter $G$. Then for any $x\in \mathbb{T}^{1}$ one may assign to the cocycle $M(x,l)$ a word $w(x)=[w_{1}, w_{2},\ldots, w_{l}]$ of length $l$, consisting of letters $B$ and $G$ such that $w_{i}=B$ if $\sigma_{\omega}^{i-1}(x)\in U_{\delta}(\mathcal{C}_{0})$ and $w_{i}=G$ otherwise. Finally, for any $x\in U_{\delta}(c_{0})$ the product $M(\sigma_{\omega}^{n_{-}}(x), n_{-}+n+n_{+}) = \prod\limits_{k=-n_{-}}^{n+n_{+}}A(\sigma_{\omega}^{k}(x))$ will be denoted by letter $H$. Note that, under resonance conditions, Lemma 6 guarantees hyperbolicity of the product $M(\sigma_{\omega}^{n_{-}}(x), n_{-}+n+n_{+})$. One may factorize the word $w(x)$ by means of the latter notation. We replace any subword $[w_{s},\ldots, w_{s+n+n_{-}+n_{+}}]$ of length $n+n_{-}+n_{+}$ by one letter $H$ if $\sigma_{\omega}^{s+n_{-}}(x)\in U_{\delta}(c_{0})$. As a result we obtain a factorized word $w^{f}(x)$ of length $j$ consisting of three letters $B, G, H$. Note that $j\to \infty$ as $l\to \infty$. Moreover, letter $B$ may appear in the word $w^{f}$ no more than two times. If letter $B$ occurs exactly two times, we define $k_{1}, k_{2}$ to be the indices such that $w^{f}_{k_{i}}=B, i=1,2$ and $k_{1}<k_{2}$. When $B$ appears less than two times, we consider three cases. If 
$w^{f}_{k}\neq B, k=1,\ldots, [j/2]$, we set $k_{1} = 0$ and $k_{2}$ such that $w^{f}_{k_{2}}=B$. If $w^{f}_{k}\neq B, k=[j/2]+1,\ldots, j$, we set $k_{2} = j+1$ and $k_{1}$ such that $w^{f}_{k_{1}}=B$. Finally, if letter $B$ does not appear in the word $w^{f}(x)$, we define $k_{1}=0, k_{2}=j+1$.
Clearly, we have
\begin{equation}\label{eqKs}
k_{1}<\tau_{0}(\delta),\quad j-k_{2}< \tau_{0}(\delta).
\end{equation}

Then one may define a truncated word $\hat w$ by the rule
\begin{equation*}
\hat w_{i} = w^{f}_{k_{1}+i},\; i=1,\ldots, k_{2}-k_{1}-1.
\end{equation*}
The truncated word corresponds to a product of matrices satisfying the conditions of Lemma 2. On the other hand, due to (\ref{eqKs})
\begin{equation*}
\frac{k_{2} - k_{1} - 1}{j}\to 1\quad  \textrm{as}\quad j\to \infty. 
\end{equation*}
Hence, for sufficiently large $j$ we may apply Lemma 3 to the cocycle $M(x,l)$. This yields the following
\begin{thms} 
Let hypotheses $(H_{1})$ - $(H_{6})$ be satisfied, $\lambda_{0}\gg 1$ and the critical set $\mathcal{C}_{0}$ consists of two points. Assume $t_{res}\in [a,b]$ is resonant of order $n$ and the time of primary collisions 
\begin{equation*}
\tau_{0} > \left(1 + \frac{2 \log \lambda_{\max}}{\log (C_{M}\lambda_{\min})}\right) n.
\end{equation*}
Then there exists $\varepsilon_{0}>0$ such that for any $\varepsilon\in (0, \varepsilon_{0})$
there exists positive $h$ such that 
\begin{equation*}
h = O\Bigl(\varepsilon^{2}\lambda_{0}^{-2n}\Bigr)
\end{equation*}
and the cocycle $M(x,n)$ is uniformly hyperbolic for all $t\in (t_{res} - h, t_{res} + h)$. 
\end{thms} 

{\bf Remark 5} It has to be noted that in Theorem 2 we do not impose any conditions on the rotation number in contrary to Theorem 1. The reason for that is the fact that those parameter values, which correspond to the uniform hyperbolicity of a cocycle, constitute an open set.

\section*{Acknowledgements}

The research was supported by RFBR grant (project No. 20-01-00451/22).


\end {document}